\documentclass{ifacconf}
\usepackage{amsmath,amssymb,amsfonts}
\usepackage{subcaption}
\usepackage{graphicx}      
\usepackage{natbib}        
\usepackage{algorithm, algcompatible}
\usepackage{algpseudocode}
\usepackage[utf8]{inputenc}
\usepackage{xcolor}

\newtheorem{assumption}{Assumption}

\newtheorem{remark}{Remark}

\newtheorem{problem}{Problem}

\begin{document}
\begin{frontmatter}

\title{Markov Parameter Identification via Chebyshev Approximation\thanksref{footnoteinfo}} 

\thanks[footnoteinfo]{This work is supported by National Natural Science Foundation of China under grant No. 62192752.}

\author[First]{Jiayun Li} 
\author[First]{Yilin Mo} 

\address[First]{Department of Automation and BNRist, Tsinghua University, Beijing, P.R.China (e-mail: lijiayun22@mails.tsinghua.edu.cn, ylmo@tsinghua.edu.cn.)}

\begin{abstract}                
This paper proposes an identification algorithm for Single Input Single Output (SISO) Linear Time-Invariant (LTI) systems. In the noise-free setting, where the first $T$ Markov parameters can be precisely estimated, all Markov parameters can be inferred by the linear combination of the known $T$ Markov parameters, of which the coefficients are obtained by solving the uniform polynomial approximation problem, and the upper bound of the asymptotic identification bias is provided. For the finite-time identification scenario, we cast the system identification problem with noisy Markov parameters into a regularized uniform approximation problem. Numerical results demonstrate that the proposed algorithm outperforms the conventional Ho-Kalman Algorithm for the finite-time identification scenario while the asymptotic bias remains negligible.
\end{abstract}


\end{frontmatter}

\section{Introduction}
Linear Time-Invariant (LTI) systems are an important class of models in many technical fields, e.g., industry~\citep{schroeck2001compensator}, automobile~\citep{krisda2012using}, robotics~\citep{gujs461494} and so on. Although conventional control methods for the LTI systems, which may not need accurate system model, are proved to be effective (e.g. the PID controller), for tasks requiring high agility and performance, precise models become a necessity. Therefore, the system identification problem has caught great attention. While classical identification research mainly focuses on the identifiability of the system, as well as the asymptotic performance of certain identification methods, e.g., the Ordinary Least Square (OLS) based method, there has recently been growing attention on the \textit{finite-time performance} of the identification algorithms, which can be further divided into two categories, depending on whether the state of the system is directly obtained. 

When the state of the system is precisely observed, the identification of system parameters using the classic Least Square (LS) method is proved to be near optimal, where several recent works are devoted to deriving concentration analysis on the identification error of the finite-time LS algorithm~\citep{shirani_faradonbeh_finite_2017, shirani_faradonbeh_finite_2018, rantzer_concentration_2018, simchowitz_learning_2018, sarkar_near_2019}. Besides the LS methods, \cite{wagenmaker_active_2020} propose an active learning method for system identification, and a theoretical upper bound on the identification error of system parameters is provided.

When the state of the system cannot be directly measured, the system identification problem becomes ``more challenging"~\citep{oymak2019non}. A great line of work aims to address the identification of \textit{a finite number of Markov parameters} from the limited sample trajectories and analyze the performance of specific identification algorithms. \cite{simchowitz_learning_2019} propose to identify finite Markov parameters using the pre-filtered least square method, and an upper bound on the identification error is further derived. Furthermore, \cite{zheng_non-asymptotic_2021} leverage the LS method to identify finite Markov parameters from the sample trajectories of both stable and unstable systems, where the performance of the LS method is further provided. In addition, another line of work approximates the \textit{transfer function} of the system by recovering Finite Impulse Response (FIR)~\citep{ljung_asymptotic_1985, ljung_asymptotic_1992, goldenshluger_nonparametric_1998} or truncated Infinite Impulse Response (IIR)~\citep{yin_maximum_2021, iannelli_experiment_2021} from the sampled data. The performance of the FIR methods are guaranteed by \cite{tu_non-asymptotic_2017} when only a limited amount of data is available. However, leveraging a finite number of Markov parameters to approximate either the input-output relationship of the system or the transfer function suffers from truncation error~\citep{tu_non-asymptotic_2017}.

On the other hand, several methods are proposed to address the above problem by identifying system parameters of the state-space model using finite Markov parameters. \cite{oymak2019non} propose the finite-time Ho-Kalman algorithm, which recovers a balanced realization of the original system by conducting SVD decomposition to the Hankel matrix formed by finite Markov parameters, and they further provide a theoretical upper bound on the identification error of system parameters. Besides the LS-based methods, subspace methods are another common approach to recover system parameters from Markov parameters, and \cite{tsiamis_finite_2019} derive an upper bound on the identification error of subspace methods. There is also a line of work that focuses on the identification performance of the classic MOESP algorithm~\citep{ikeda_estimation_2015, chiuso_ill-conditioning_2004}, which is a well-known subspace method. 

However, our initial work~\citep{li2022fundamental} provides theoretical analysis on the sample complexity of the Ho-Kalman algorithm with respect to system dimension, and prove that the algorithm is ill-conditioned for high-dimensional systems. Besides, the ill-conditionedness of subspace methods is also observed both in theory and practice, e.g., see~\cite{chiuso_ill-conditioning_2004,hachicha2014n4sid}. As a result, although the above algorithms are proved to be asymptotically unbiased, the finite-time performance may be poor, especially for identification problems of high-dimensional systems. Numerical results further demonstrate the finite-time performance of the Ho-Kalman algorithm in Section~\ref{sec:numerical_result}.

Besides the specific identification algorithms, our initial works also prove that the ill-conditionedness is actually rooted in the identification problem itself, and that the identification problem is ill-conditioned for high-dimensional systems using \textit{any unbiased identification algorithms}~\citep{sun2022fundamental,sun2022identification,li2022fundamental}. 

This result motivates us to avoid identifying state-space parameters, and instead derive a biased Markov parameter identification algorithm. We circumvent ill-conditionedness by minimizing the mean squared identification error of Markov parameters using Chebyshev approximation. Taking notice of the fact that infinite Markov parameters contain full information about the input-output relationship of the original system, we choose to recover all Markov parameters using the proposed algorithm. Furthermore, an upper bound on the asymptotic identification bias is provided, which indicates that the bias is negligible compared to the impact of noise in the original data.


The main contribution of this paper is as follows:
\begin{itemize}
   \item We formulate the Markov parameter identification as a uniform polynomial approximation problem, assuming the first $T$ Markov parameters can be precisely observed.
   \item When only noisy Markov parameters from finite sample trajectories are available, we further cast the identification problem into a regularized polynomial approximation problem.
   \item An upper bound on the asymptotic identification bias of the proposed algorithm is provided, which is negligible compared to the impact of noise in the original data.
\end{itemize}

The rest of the paper is organized as follows: In Section~\ref{sec:problem_formulation}, we formulate the identification problem in the SISO case. In addition, Section~\ref{sec:noise_free_case} illustrates our algorithm in the noise-free setting, and a theoretical upper bound on the identification bias is provided. Section~\ref{sec:noise_case} analyzes the mean squared identification error in the noisy scenario, by which we formulate the finite-time identification problem as regularized uniform approximation. Numerical results are depicted in Section~\ref{sec:numerical_result}, which verifies the finite-time performance of our algorithm, and Section~\ref{sec:conclusion} concludes the paper.

\section{Problem Formulation}\label{sec:problem_formulation}

Consider a Single Input Single Output (SISO) discrete-time LTI system with the following state-space realization:
\begin{equation}\label{eq:linear_system}
	\begin{aligned}
		x_{t+1} &= Ax_t + Bu_t+w_t, \\
		y_t &= Cx_t + v_t,
	\end{aligned}
\end{equation}
where $ x_t \in \mathbb{R}^n,  u_t \in \mathbb{R}, y_t \in \mathbb{R}  $ are the system state, input and output at time $ t $ respectively, and $n$ is the dimension of the system. $\{w_t\}$ and $\{v_t\}$ are independent and identically distributed (i.i.d.) Gaussian noise with zero mean and covariance $Q$, $R$, (with $Q, R\geq 0$), respectively. $A, B$ and $C$ are system parameters with proper dimensions.

The overarching goal of this paper is to identify all Markov parameters $\{H_t\}_{t=1}^\infty$, which contain full information about the input-output relationship of the original system, from the first $T$ exact Markov parameters $\{H_t\}_{t=1}^{T}$ in Section~\ref{sec:noise_free_case} or from the estimated Markov parameters $\{\tilde H_t\}_{t=1}^T$ in Section~\ref{sec:noise_case}.

We further make the following assumptions:
\begin{assumption}\label{ass:main}
      System matrix $A$ is real diagonalizable, and all eigenvalues of $A$ lie on the real interval $[-\rho, \rho]$, where $\rho>0$ is known.
\end{assumption}

Since similarity transformation preserves the input-output relationship of the original system, for convenience of notation, we assume that the system takes the diagonal canonical form, i.e., $A=\text{diag}\{\lambda_1, \lambda_2, \cdots, \lambda_n\}$, where $\lambda_i, i=1, \cdots, n$ are the eigenvalues of system matrix $A$. $B$ is an all one column vector, i.e., $B=\begin{bmatrix} 1 & 1 & \cdots & 1\end{bmatrix}^\top$. $C=\begin{bmatrix}c_1 & c_2 & \cdots & c_n\end{bmatrix}$, and we further define the modified system energy as 
\begin{equation}\label{eq:cm_def}
   C_m\triangleq \|C\|_1.
\end{equation}
Note that though the roots $\lambda_i$ are denoted distinctly, the eigenvalues with different indices can have the same value.
\begin{remark}
   The identification of state-space parameters $A, B$ and $C$ and the estimation of all Markov parameters are two independent schemes that are able to recover full input-output relationship of the original system. Previous identification algorithms mainly consider the estimation of state-space parameters, while we choose to directly identify all Markov parameters of the LTI system.
\end{remark}
\begin{remark}
   For simplicity and convenience of notation, we only discuss the identification of SISO systems in this paper. However, it is worth noting that the proposed algorithm is directly applicable to Multiple-Input Multiple-Output (MIMO) systems without significant modification.
\end{remark}
\begin{remark}
   The identification of systems with complex poles is left for future research.
\end{remark}

\section{Identification algorithm in the noise-free case}\label{sec:noise_free_case}
This section considers the identification problem in the noise-free case, where we make the following assumption:
\begin{assumption}
   The first $T$ Markov parameters $\{H_t\}_{t=1}^T$ are perfectly known.
\end{assumption}

The proposed algorithm aims to identify $H_{k}$ through the \textit{linear combination} of the first $T$ Markov parameters $\{H_t\}_{t=1}^T$, where $k>T$:
\begin{equation}\label{eq:linear_combination}
   \hat H_{k}=\sum_{t=1}^T\alpha_{t-1}H_t,
\end{equation}
such that the mean squared identification error is minimized:
\begin{equation}\label{eq:initial_prob}
   \epsilon\triangleq |H_{k}-\hat H_{k}|^2=\left|H_k-\sum_{t=1}^T\alpha_{t-1}H_t\right|^2.
\end{equation}

According to Cayley-Hamilton Theorem, for each $k$, $H_k$ can be expressed as a linear combination of the first $n$ Markov parameters. Thus, $\epsilon=0$ can be reached by choosing system-specific coefficients $\alpha_t$ when system parameters are known, which is not satisfied by the identification problem. On the other hand, noticing that the $k$-th Markov parameter of system~\eqref{eq:linear_system} has the following relationship with state-space parameters:
\begin{equation}\label{eq:Markov_parameter_form}
   H_{k}=CA^{k-1}B=\sum_{j=1}^n c_j\lambda_j^{k-1},
\end{equation}
we instead formulate the minimization of $\epsilon$ in~\eqref{eq:initial_prob} as an approximation problem below:

\begin{problem}\label{prob:poly_approx_prob}
   \begin{equation}\label{eq:unregularized_uniform_approx}
      \boldsymbol{\alpha}^*=\text{argmin}_{\{\alpha_t\}_{t=0}^{T-1}} \left\|\lambda^{k-1}-\sum_{t=0}^{T-1} \alpha_{t}\lambda^t\right\|_{\infty, \rho},
   \end{equation}
   in which $\|\cdot\|_{\infty, \rho}$ denotes the uniform norm of functions, i.e., $\|f\|_{\infty, \rho}=\sup\{|f(\lambda)|: \lambda\in[-\rho, \rho]\}$, $\boldsymbol{\alpha} ^*\triangleq \begin{bmatrix} \alpha^*_0 & \alpha^*_1 & \cdots & \alpha^*_{T-1}\end{bmatrix}^\top$ is the optimal solution.
\end{problem}
Therefore, $\boldsymbol{\alpha}^*$ becomes a universal solution without the dependence on system parameters $A, B, C$ and $n$, and is only associated with $T, \rho$ and $k$. 

Furthermore, the estimation of the $k$-th Markov parameter can be calculated as:
\begin{equation}
   \begin{aligned}
      \hat H_k=\sum_{t=1}^T \alpha_{t-1}^* H_t.
   \end{aligned}
\end{equation}



\subsection{Performance Analysis}\label{subsec:noise_free_performance}
In this subsection, we analyze the performance of the above algorithm, where the proof of the following theorem is provided in Appendix~\ref{append:asymp_error}.
\begin{thm}\label{thm:asymp_error}
   The mean squared identification error of the $k$-th Markov parameter of system~\eqref{eq:linear_system} using the first $T$ Markov parameters $\{H_t\}_{t=1}^T$ with the algorithm in~\eqref{eq:unregularized_uniform_approx} has the following upper bound, $\forall k>T$:
   \begin{equation}\label{eq:asymp_error}
      \begin{aligned}
         & \left|H_{k}-\hat H_{k}\right|^2\leq C_m^2\rho^{2k-2}\min\left(4\exp\left(-\frac{2(T-1)^2}{k-1}\right), \frac{1}{4}\right).
      \end{aligned}
   \end{equation}
   where $\hat H_{k}=\sum_{t=1}^T \alpha_{t-1}^*H_t$ is the estimation of $H_{k}$, the coefficients $\boldsymbol{\alpha}^*$ are obtained by~\eqref{eq:unregularized_uniform_approx} and $C_m$ is the modified system energy defined in~\eqref{eq:cm_def}.
\end{thm}

\begin{remark}\label{remark:noise_free_performance}
   When the system is strictly stable and $\rho<1$, the \textit{finite} $\mathcal{L}_2$ identification error of all Markov parameters can be obtained by putting all estimation of Markov parameters to be $0$, i.e., $\hat H_k=0, \forall k>T$ (which shall be referred to as \textit{the truncation method} in the following): 
   \begin{equation}
      \begin{aligned}
         \sqrt{\sum_{k=1}^\infty |H_k-\hat H_k|^2}&=\sqrt{\sum_{k=T+1}^\infty \left|\sum_{j=1}^n c_j\lambda_j^{k-1}\right|^2}\\
         &\leq \sqrt{\sum_{k=T}^\infty C_m^2\rho^{2k}}=\frac{C_m\rho^T}{\sqrt{1-\rho^2}}.
      \end{aligned}
   \end{equation}
   Furthermore, Theorem~\ref{thm:asymp_error} reveals that the upper bound on the identification error of each Markov parameter is better than the truncation method:
   \begin{equation}
      \begin{aligned}
         |H_k-\hat H_k|^2&\leq C_m^2\rho^{2k-2}\min \left(4\exp\left(-\frac{2(T-1)^2}{k-1}\right), \frac{1}{4}\right) \\
         & \leq C_m^2\rho^{2k-2}, \forall k > T,
      \end{aligned}
   \end{equation}
   which indicates that the $\mathcal{L}_2$ identification error of the proposed algorithm is also finite and is strictly better than $\frac{C_m\rho^T}{\sqrt{1-\rho^2}}$. Moreover, the exponential term $\exp\left(-\frac{2(T-1)^2}{k-1}\right)$ in~\eqref{eq:asymp_error} guarantees that the upper bound on the identification error of the first few unknown Markov parameters, which dominates the $\mathcal{L}_2$ error of the truncation method, is drastically reduced by our algorithm.

   On the other hand, when the system is unstable, since the value of Markov parameter $H_k$ explodes to infinity with the increase of the index $k$, only a finite number of Markov parameters are worth identifying, where the identification error is also bounded by Theorem~\ref{thm:asymp_error}. 
\end{remark}
\begin{remark}
   It is worth noting that the identification bias can be efficiently reduced by increasing the number of known Markov parameters $T$, which is further illustrated considering the following three cases:
   \begin{itemize}
      \item For a fixed index $k$, the identification error of $H_{k}$ attenuates superlinearly w.r.t. the number of known Markov parameters $T$:
      \begin{equation}
         |H_{k}-\hat H_{k}|^2\leq 4C_m^2\rho^{2k-2}\exp\left(-\frac{2(T-1)^2}{k-1}\right).
      \end{equation}
      \item When the system is strictly stable and $\rho<1$, the identification bias has the following supremum among all Markov parameters:
      \begin{equation}\label{eq:mse_error_sup}
         \sup_k |H_k-\hat H_k|^2\leq 4C_m^2\exp\left(-4(T-1)\sqrt{\log 1/\rho}\right),
      \end{equation}
      which attenuates exponentially w.r.t. $T$.
      
      \item When the system is strictly stable and $\rho<1$, in order to obtain identification accuracy $\delta^2$, by~\eqref{eq:mse_error_sup}, the number of known Markov parameters required satisfies
      \begin{equation}
         \begin{aligned}
            T\sim\mathcal{O}\left(\frac{\log(2C_m)-\log\delta}{\sqrt{\log(1/\rho)}}\right),
         \end{aligned}
      \end{equation}
      and when $\rho$ is close to $1$,
      \begin{equation}
         T\sim\mathcal{O}\left(\frac{\log(2C_m)-\log\delta}{\sqrt{1-\rho}}\right).
      \end{equation}
      The above result indicates that the number of known Markov parameters $T$ required for a specific system increases at an efficient speed of $T\sim\mathcal{O}\left(\log\left(\frac{1}{\delta}\right)\right).$
   \end{itemize}
\end{remark}


\section{Identification in the noisy scenario}\label{sec:noise_case}
In this section, we address the identification problem where only a finite amount of data is available, and the estimate of Markov parameters $\{\tilde H_t\}_{t=1}^T$ is contaminated by noise. We make the following assumption on the obtained data:
\begin{assumption}
   The sample trajectories are collected episodically from the system, where each episode is reset to steady state and is independent from each other.
\end{assumption}

\subsection{Mean Squared Identification Error}\label{subsec:noisy_error}
We first illustrate one method to recursively identify $\tilde H_t$, where $t\leq T$, using the sample trajectories obtained from the episodic experiments.

In each episode of the experiment, we inject the \textit{unit impulse function} to the system, i.e.,
\begin{equation}
   u_t=\left\{\begin{aligned}
      1,\quad &t=0, \\
      0, \quad&t\geq 1,
   \end{aligned}\right.
\end{equation}
and we denote the $t$-th output in the $\ell$-th episode as $y_t^{(\ell)}, t=1, \cdots, T$. The estimation of the $t$-th Markov parameter can be updated recursively by: 
\begin{equation}
   \tilde H_t^{(\ell)}=\tilde H_t^{(\ell-1)}+\frac{1}{\ell}\left(y_t^{(\ell)}-\tilde H_t^{(\ell-1)}\right), \tilde H_t^{(0)}=0.
\end{equation}

The following theorem quantifies the mean squared identification error of $H_k$ using the linear combination of the first $T$ Markov parameters  $\{\tilde H_t^{(N)}\}_{t=1}^T$ estimated from $N$ sample trajectories, where the proof is reported in Appendix~\ref{append:markov_error}.
\begin{thm}\label{thm:markov_approx_error}
   Let $\hat H_k=\sum_{t=1}^T \alpha_{t-1}\tilde H_t^{(N)}$ be the approximation of the $k$-th Markov parameter $H_k (k>T)$. Then, the mean squared approximation error of the $k$-th Markov parameter has the following upper bound:
   \begin{equation}\label{eq:markov_approx_error}
      \mathbb{E}\left(H_k-\hat H_k\right)^2\leq C_m^2\left\|\lambda^{k-1}-\sum_{t=0}^{T-1} \alpha_{t}\lambda^t\right\|_{\infty, \rho}^2+\frac{\Sigma}{N}\|\boldsymbol{\alpha}\|_1^2,
   \end{equation}
   where $\Sigma$ is the variance of the estimation of Markov parameters: 
   \begin{equation}\label{eq:markov_para_variance}
      \Sigma=\mathbb{E}\left(H_t-y_t^{(\ell)}\right)^2=CPC^\top+R, \text{ where }P=APA^\top+Q,
   \end{equation}
   $\forall \ell=1, \cdots, N, t=1, \cdots, T$, $A, C, Q$ and $R$ are system parameters defined in Section~\ref{sec:problem_formulation}.
\end{thm}

\begin{remark}
   The assumption that each experiment is reset to steady state is only for simplicity of the analysis in Theorem~\ref{thm:markov_approx_error}. Generally, the proposed algorithm is able to work without this assumption.
\end{remark}
\begin{remark}
   Although $\Sigma$ cannot be calculated by~\eqref{eq:markov_para_variance} since the system parameters are unknown, it can be estimated from data using the following equation:
   \begin{equation}
      \hat\Sigma=\frac{1}{TN-1}\sum_{t=1}^T\sum_{\ell=1}^N\left(y_t^{(\ell)}-\bar y_t\right)^2,
   \end{equation}
   where $\bar y_t=\frac{1}{N}\sum_{\ell=1}^N y_t^{(\ell)}$.
\end{remark}
\begin{remark}
   Note that the proposed estimation scheme regarding the first $T$ Markov parameters can be replaced by other Markov parameter estimators.
\end{remark}


\vspace{-0.2cm}
\subsection{Identification Algorithm}\label{subsec:noisy_alg}
We now propose our regularized identification algorithm in the noisy scenario, which aims to minimize the mean squared identification error of the Markov parameter $H_k$:
\begin{problem}\label{prob:final_opt_prob}
   \begin{equation}\label{eq:regularized_opt_prob}
      \boldsymbol{\alpha}^*=\text{argmin}_{\{\alpha_t\}_{t=0}^{T-1}} \left\|\lambda^{k-1}-\sum_{t=0}^{T-1}\alpha_t\lambda^t\right\|_{\infty, \rho}^2+\gamma\|\boldsymbol{\alpha}\|_1^2,
   \end{equation}
\end{problem}
in which $\gamma$ is a parameter indicating the relative importance between the above two terms, and shall be chosen as $\gamma=\frac{\hat\Sigma}{C_m^2N}$ according to~\eqref{eq:markov_approx_error}.

Suppose $N$ sample trajectories are leveraged to estimate the first $T$ Markov parameters $\{H_t\}_{t=1}^T$, the inference of $H_k$ can be calculated as:
\begin{equation}
   \hat H_k=\left\{ \begin{aligned}
      &\tilde H_k^{(N)}, &k\leq T; \\
      &\sum_{t=1}^T \alpha_{t-1}^*\tilde H_t^{(N)}, &k>T.
   \end{aligned}\right.
\end{equation}

Please refer to Algorithm~\ref{alg:alg} for a detailed illustration of the proposed algorithm.

\begin{algorithm}
   \caption{Complete procedure of our regularized identification algorithm proposed in Section~\ref{subsec:noisy_alg}.}\label{alg:alg}
   \begin{algorithmic}
      \Require The index of Markov parameter to be identified $k\geq 1$, maximal steps of each episode $T$, total iteration epochs $N$, modified system energy $C_m$
      \Ensure Identified $k$-th Markov parameter $\hat H_k$
      \State Initialization $\tilde H_1\leftarrow 0, \cdots, \tilde H_T\leftarrow 0$, $\ell\leftarrow 1$
      \While{$\ell\leq N$}\Comment{Identification of the first $T$ Markov parameters}
         \State $\tilde H_t\leftarrow\tilde H_t+\frac{1}{\ell}\left(y_t^{(\ell)}-\tilde H_t\right), t=1, \cdots, T$
         \State $\ell\leftarrow \ell+1$
      \EndWhile
      \If{$k>T$}\Comment{Identification of $H_k$ via polynomial approximation}
         \State $\bar y_t\leftarrow\frac{1}{N}\sum_{\ell=1}^N y_t^{(\ell)}, t=1, \cdots, T$
         \State $\hat\Sigma\leftarrow\frac{1}{TN-1}\sum_{t=1}^T\sum_{\ell=1}^N\left(y_t^{(\ell)}-\bar y_t\right)^2$
         \State $\gamma\leftarrow\frac{\hat\Sigma}{C_m^2N}$
         \State $\boldsymbol{\alpha}^*\leftarrow\text{argmin}_{\{\alpha_t\}_{t=0}^{T-1}} \left\|\lambda^{k-1}-\sum_{t=0}^{T-1}\alpha_t\lambda^t\right\|_{\infty, \rho}^2+\gamma\|\boldsymbol{\alpha}\|_1^2$
         \State $\hat H_k\leftarrow \sum_{t=1}^T \alpha_{t-1}^*\tilde H_t$
      \Else
         \State $\hat H_k\leftarrow \tilde H_k$
      \EndIf
   \end{algorithmic}
\end{algorithm}

\section{Numerical Results}\label{sec:numerical_result}
This section compares the finite-time performance of the unbiased Ho-Kalman algorithm~\citep{oymak2019non}, the truncation method introduced in Remark~\ref{remark:noise_free_performance} and the proposed algorithm using numerical experiments.

\begin{figure}
   \centering
   \begin{subfigure}{\columnwidth}
      \centering
      \includegraphics[width=\columnwidth]{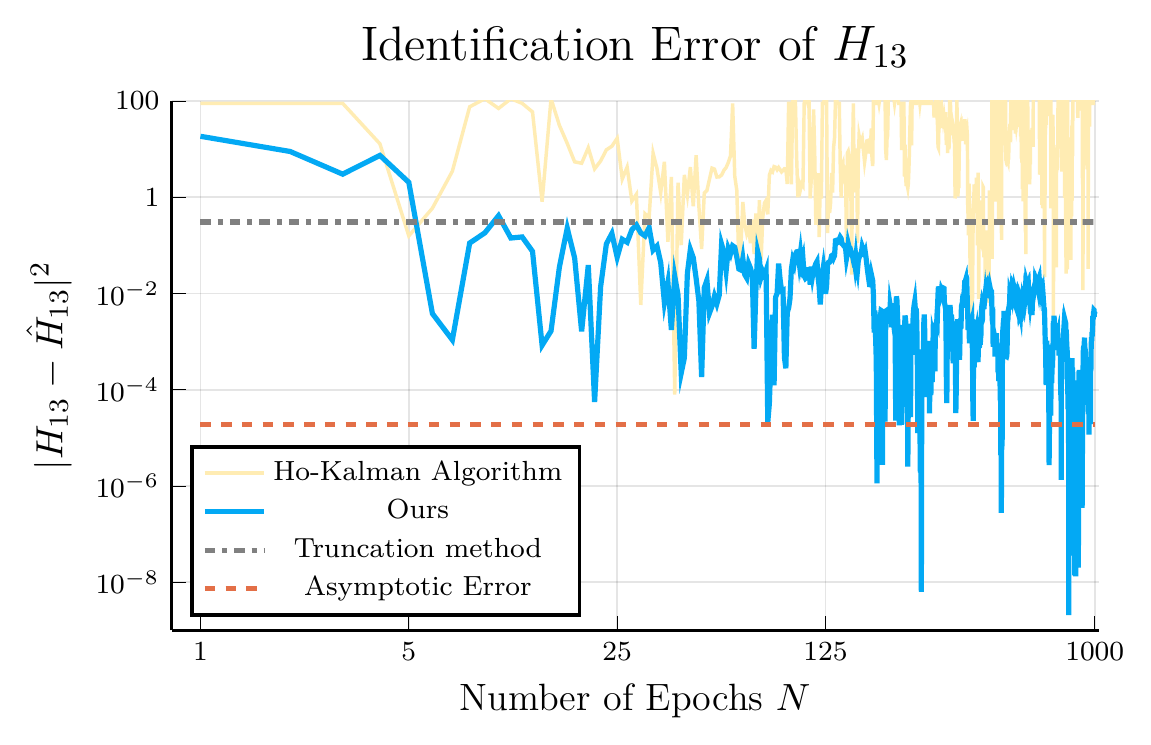}
   \end{subfigure}
   \begin{subfigure}{\columnwidth}
      \centering
      \includegraphics[width=\columnwidth]{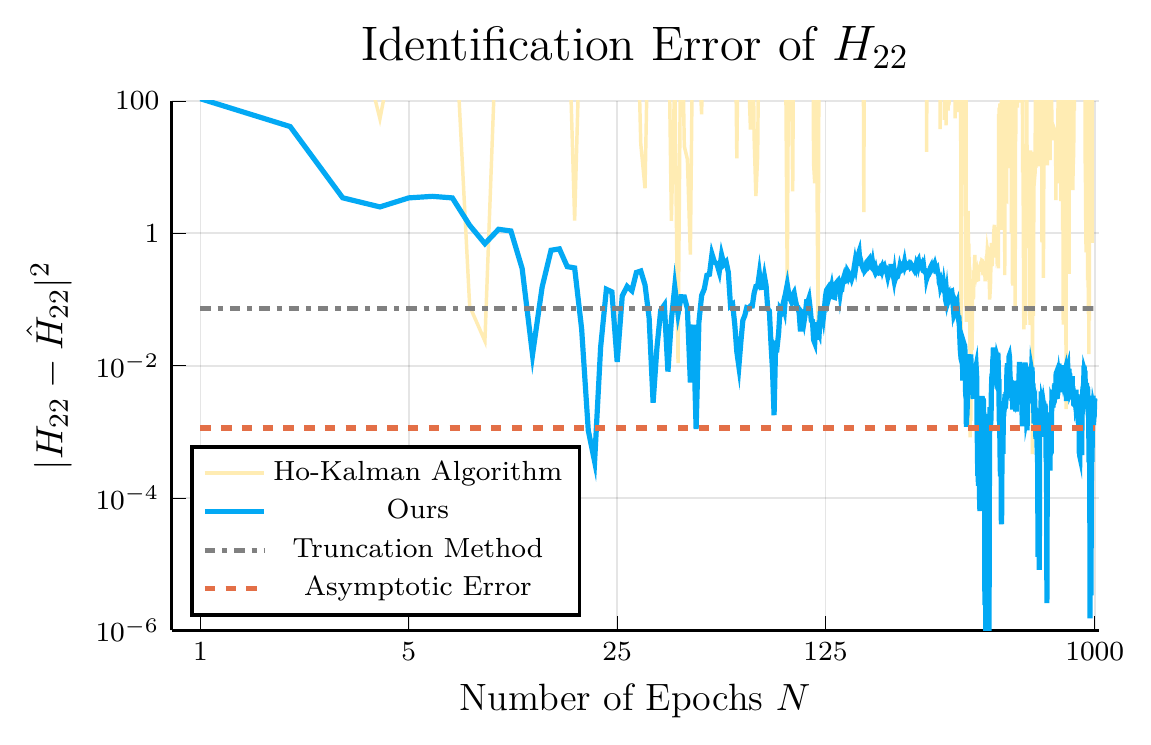}
   \end{subfigure}
   \caption{Identification error of $H_{13}$ and $H_{22}$ using $\{\tilde H_{i}\}_{i=1}^{12}$ by the Ho-Kalman algorithm, the proposed algorithm and the truncation method respectively, which is depicted in a log-log plot. The yellow, blue, and gray lines depict the identification error of the Ho-Kalman algorithm, the proposed algorithm, and the truncation method, respectively. Furthermore, the orange line denotes the asymptotic identification bias of the proposed algorithm when $N$ tends to infinity, and the weighting parameter in~\eqref{eq:regularized_opt_prob} tends to $0$.}
   \label{fig:exp1_result}
\end{figure}

We consider the following $6$-dimensional SISO system with $\rho=0.95$:
\begin{equation}\label{eq:exp_1_sys}
   \begin{aligned}
      A& =\text{diag}(0.94, 0.75, -0.75, -0.69, 0.46, 0.42), \\
      B& =\begin{bmatrix}
         1 & 1 & 1 & 1 & 1 & 1
      \end{bmatrix}^\top, C=\begin{bmatrix}
         1 & 1 & 1 & 1 & 1 & 1
      \end{bmatrix}.
   \end{aligned}
\end{equation}
We iteratively estimate the target Markov parameter after each episode, and let the weighting coefficient $\gamma$ in~\eqref{eq:regularized_opt_prob} attenuate at a rate of $\mathcal{O}\left(\frac{1}{N}\right)$, where $N$ is the number of epochs, to illustrate the finite-time performance of our algorithm with the growth of the amount of data sampled. The identification result of $H_{13}$ and $H_{22}$ using $\{\tilde H_i\}_{i=1}^{12}$ by our algorithm, the Ho-Kalman algorithm, and the truncation method are compared in Figure~\ref{fig:exp1_result}.

As can be observed from Figure~\ref{fig:exp1_result}, although the Ho-Kalman algorithm is proved to be asymptotically unbiased, the finite-time identification result oscillates to at least $10$ times larger than the true value of $H_{13}$, and shows little decay during the first $1000$ episodes. On the contrary, our algorithm smoothly identifies the Markov parameter and quickly converges to the asymptotic value. Our algorithm also outperforms the naive truncation method, which reveals the effectiveness of the proposed optimization problem.

Similar results are obtained from the identification process of $H_{22}$ plotted in Figure~\ref{fig:exp1_result}, where our algorithm significantly outperforms the Ho-Kalman algorithm. In contrast, the bias of our algorithm remains negligible compared to the impact of noise in the original data. Note that the asymptotic bias can be further reduced by slightly increasing the number of known Markov parameters $T$, where the impact of $T$ on the identification bias is illustrated in Section~\ref{subsec:noise_free_performance}.

\begin{figure}
   \centering
   \includegraphics[width=\columnwidth]{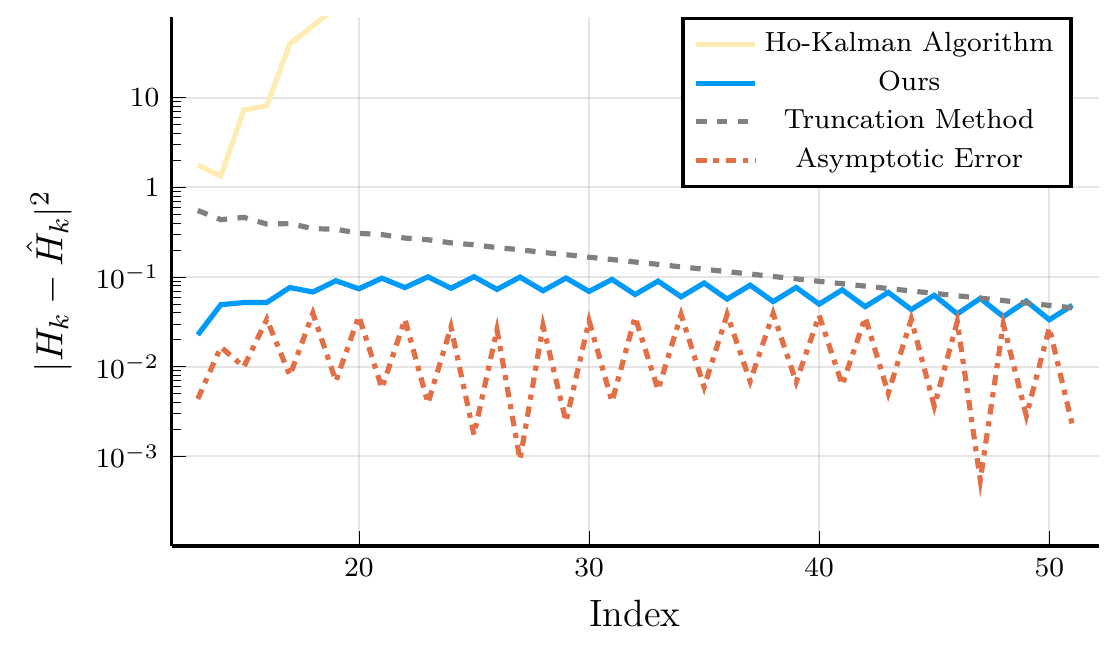}
   \caption{Identification error of the Ho-Kalman algorithm, the proposed algorithm and the truncation method w.r.t. the index of Markov parameters using the first $12$ Markov parameters from $1000$ sample trajectories of system~\eqref{eq:exp_1_sys} in a single experiment, which is depicted in a log plot. The yellow, blue, and gray lines depict the identification error of the Ho-Kalman algorithm, the proposed algorithm, and the truncation method, respectively. The orange line denotes the asymptotic identification error of the proposed algorithm when $N$ tends to infinity, and the weighting parameter in~\eqref{eq:regularized_opt_prob} tends to $0$.}
   \label{fig:exp3_result}
\end{figure}

Finally, the identification error of $\{H_i\}_{i=13}^{50}$ using the first $12$ Markov parameters estimated from $1000$ sample trajectories of system~\eqref{eq:exp_1_sys} \textit{in a single experiment} are shown in Figure~\ref{fig:exp3_result}. The result reveals that the identification error of the Ho-Kalman algorithm explodes exponentially w.r.t. the index of Markov parameters due to the fact that some poles of the identified $A$ matrix are unstable. On the other hand, the identification result of our algorithm outperforms both the Ho-Kalman algorithm and the naive truncation method, while the bias of our algorithm remains imperceptible.

\section{Conclusion}\label{sec:conclusion}
This paper proposes a Markov parameter identification algorithm for SISO LTI systems from the first $T$ Markov parameters, which is suitable for both the noise-free and the noisy scenarios. We further provide the upper bound on the asymptotic identification bias when the first $T$ Markov parameters are precisely estimated. In addition, the mean squared identification error in the noisy scenario is derived, by which we formulate the finite-time identification problem as a regularized uniform polynomial approximation problem. The identification algorithm for systems with complex roots is left for future research.


\bibliography{ref}             

\begin{thebibliography}{26}
\providecommand{\natexlab}[1]{#1}
\providecommand{\url}[1]{\texttt{#1}}
\providecommand{\urlprefix}{URL }
\expandafter\ifx\csname urlstyle\endcsname\relax
  \providecommand{\doi}[1]{doi:\discretionary{}{}{}#1}\else
  \providecommand{\doi}{doi:\discretionary{}{}{}\begingroup
  \urlstyle{rm}\Url}\fi

\bibitem[{Chiuso and Picci(2004)}]{chiuso_ill-conditioning_2004}
Chiuso, A. and Picci, G. (2004).
\newblock On the ill-conditioning of subspace identification with inputs.
\newblock \emph{Automatica}, 40(4), 575--589.
\newblock \doi{10.1016/j.automatica.2003.11.009}.
\newblock
  \urlprefix\url{https://linkinghub.elsevier.com/retrieve/pii/\\S0005109803003674}.

\bibitem[{Faradonbeh et~al.(2017)Faradonbeh, Tewari, and
  Michailidis}]{shirani_faradonbeh_finite_2017}
Faradonbeh, M.K.S., Tewari, A., and Michailidis, G. (2017).
\newblock Finite time analysis of optimal adaptive policies for
  linear-quadratic systems.
\newblock \emph{arXiv preprint arXiv:1711.07230}.

\bibitem[{Goldenshluger(1998)}]{goldenshluger_nonparametric_1998}
Goldenshluger, A. (1998).
\newblock Nonparametric estimation of transfer functions: rates of convergence
  and adaptation.
\newblock \emph{IEEE Transactions on Information Theory}, 44(2), 644--658.
\newblock \doi{10.1109/18.661510}.
\newblock Conference Name: IEEE Transactions on Information Theory.

\bibitem[{Hachicha et~al.(2014)Hachicha, Kharrat, and
  Chaari}]{hachicha2014n4sid}
Hachicha, S., Kharrat, M., and Chaari, A. (2014).
\newblock N4sid and moesp algorithms to highlight the ill-conditioning into
  subspace identification.
\newblock \emph{International Journal of Automation and Computing}, 11(1),
  30--38.

\bibitem[{Iannelli et~al.(2021)Iannelli, Yin, and
  Smith}]{iannelli_experiment_2021}
Iannelli, A., Yin, M., and Smith, R.S. (2021).
\newblock Experiment design for impulse response identification with signal
  matrix models.
\newblock \urlprefix\url{http://arxiv.org/abs/2012.08126}.
\newblock ArXiv:2012.08126 [cs, eess].

\bibitem[{Ikeda and Oku(2015)}]{ikeda_estimation_2015}
Ikeda, K. and Oku, H. (2015).
\newblock Estimation error analysis of system matrices in some subspace
  identification methods.
\newblock In \emph{2015 10th {Asian} {Control} {Conference} ({ASCC})}, 1--6.
\newblock \doi{10.1109/ASCC.2015.7244557}.

\bibitem[{Kritayakirana and Gerdes(2012)}]{krisda2012using}
Kritayakirana, K. and Gerdes, J.C. (2012).
\newblock Using the centre of percussion to design a steering controller for an
  autonomous race car.
\newblock \emph{Vehicle System Dynamics}, 50(sup1), 33--51.
\newblock \doi{10.1080/00423114.2012.672842}.
\newblock \urlprefix\url{https://doi.org/10.1080/00423114.2012.672842}.

\bibitem[{Li et~al.(2022)Li, Sun, and Mo}]{li2022fundamental}
Li, J., Sun, S., and Mo, Y. (2022).
\newblock Fundamental limit on siso system identification.
\newblock In \emph{2022 IEEE 61st Conference on Decision and Control (CDC)},
  856--861.
\newblock \doi{10.1109/CDC51059.2022.9993203}.

\bibitem[{Ljung and Wahlberg(1992)}]{ljung_asymptotic_1992}
Ljung, L. and Wahlberg, B. (1992).
\newblock Asymptotic {Properties} of the {Least}-{Squares} {Method} for
  {Estimating} {Transfer} {Functions} and {Disturbance} {Spectra}.
\newblock \emph{Advances in Applied Probability}, 24(2), 412--440.
\newblock \doi{10.2307/1427698}.
\newblock \urlprefix\url{https://www.jstor.org/stable/1427698}.
\newblock Publisher: Applied Probability Trust.

\bibitem[{Ljung and Yuan(1985)}]{ljung_asymptotic_1985}
Ljung, L. and Yuan, Z.D. (1985).
\newblock Asymptotic properties of black-box identification of transfer
  functions.
\newblock \emph{IEEE Transactions on Automatic Control}, 30(6), 514--530.
\newblock \doi{10.1109/TAC.1985.1103995}.
\newblock Conference Name: IEEE Transactions on Automatic Control.

\bibitem[{Oymak and Ozay(2019)}]{oymak2019non}
Oymak, S. and Ozay, N. (2019).
\newblock Non-asymptotic identification of lti systems from a single
  trajectory.
\newblock In \emph{2019 American control conference (ACC)}, 5655--5661. IEEE.

\bibitem[{Rantzer(2018)}]{rantzer_concentration_2018}
Rantzer, A. (2018).
\newblock Concentration {Bounds} for {Single} {Parameter} {Adaptive} {Control}.
\newblock In \emph{2018 {Annual} {American} {Control} {Conference} ({ACC})},
  1862--1866.
\newblock \doi{10.23919/ACC.2018.8431891}.
\newblock ISSN: 2378-5861.

\bibitem[{Saibaba(2021)}]{saibaba2021approximating}
Saibaba, A.K. (2021).
\newblock Approximating monomials using chebyshev polynomials.
\newblock \emph{arXiv preprint arXiv:2101.06818}.

\bibitem[{Sarkar and Rakhlin(2019)}]{sarkar_near_2019}
Sarkar, T. and Rakhlin, A. (2019).
\newblock Near optimal finite time identification of arbitrary linear dynamical
  systems.
\newblock \emph{arXiv:1812.01251 [cs]}.
\newblock \urlprefix\url{http://arxiv.org/abs/1812.01251}.
\newblock ArXiv: 1812.01251.

\bibitem[{Schroeck et~al.(2001)Schroeck, Messner, and
  McNab}]{schroeck2001compensator}
Schroeck, S.J., Messner, W.C., and McNab, R.J. (2001).
\newblock On compensator design for linear time-invariant dual-input
  single-output systems.
\newblock \emph{IEEE/ASME Transactions on mechatronics}, 6(1), 50--57.

\bibitem[{Shirani~Faradonbeh et~al.(2018)Shirani~Faradonbeh, Tewari, and
  Michailidis}]{shirani_faradonbeh_finite_2018}
Shirani~Faradonbeh, M.K., Tewari, A., and Michailidis, G. (2018).
\newblock Finite time identification in unstable linear systems.
\newblock \emph{Automatica}, 96, 342--353.
\newblock \doi{10.1016/j.automatica.2018.07.008}.
\newblock
  \urlprefix\url{https://linkinghub.elsevier.com/retrieve/pii/\\S0005109818303546}.

\bibitem[{Simchowitz et~al.(2019)Simchowitz, Boczar, and
  Recht}]{simchowitz_learning_2019}
Simchowitz, M., Boczar, R., and Recht, B. (2019).
\newblock Learning {Linear} {Dynamical} {Systems} with {Semi}-{Parametric}
  {Least} {Squares}.
\newblock In \emph{Proceedings of the {Thirty}-{Second} {Conference} on
  {Learning} {Theory}}, 2714--2802. PMLR.
\newblock
  \urlprefix\url{https://proceedings.mlr.press/v99/simchowitz\\19a.html}.
\newblock ISSN: 2640-3498.

\bibitem[{Simchowitz et~al.(2018)Simchowitz, Mania, Tu, Jordan, and
  Recht}]{simchowitz_learning_2018}
Simchowitz, M., Mania, H., Tu, S., Jordan, M.I., and Recht, B. (2018).
\newblock Learning {Without} {Mixing}: {Towards} {A} {Sharp} {Analysis} of
  {Linear} {System} {Identification}.
\newblock \emph{arXiv:1802.08334 [cs, math, stat]}.
\newblock \urlprefix\url{http://arxiv.org/abs/1802.08334}.
\newblock ArXiv: 1802.08334.

\bibitem[{Sun and Mo(2022)}]{sun2022identification}
Sun, S. and Mo, Y. (2022).
\newblock Fundamental identification limit on diagonal canonical form for siso
  systems.
\newblock In \emph{2022 IEEE 17th International Conference on Control \&
  Automation (ICCA)}, 728--733.
\newblock \doi{10.1109/ICCA54724.2022.9831907}.

\bibitem[{Sun et~al.(2022)Sun, Mo, and You}]{sun2022fundamental}
Sun, S., Mo, Y., and You, K. (2022).
\newblock Fundamental identification limit of single-input and single-output
  linear time-invariant systems.
\newblock In \emph{2022 13th Asian Control Conference (ASCC)}, 2157--2162.
\newblock \doi{10.23919/ASCC56756.2022.9828137}.

\bibitem[{Tsiamis and Pappas(2019)}]{tsiamis_finite_2019}
Tsiamis, A. and Pappas, G.J. (2019).
\newblock Finite {Sample} {Analysis} of {Stochastic} {System} {Identification}.
\newblock In \emph{2019 {IEEE} 58th {Conference} on {Decision} and {Control}
  ({CDC})}, 3648--3654. IEEE, Nice, France.
\newblock \doi{10.1109/CDC40024.2019.9029499}.
\newblock \urlprefix\url{https://ieeexplore.ieee.org/document/9029499/}.

\bibitem[{Tu et~al.(2017)Tu, Boczar, Packard, and
  Recht}]{tu_non-asymptotic_2017}
Tu, S., Boczar, R., Packard, A., and Recht, B. (2017).
\newblock Non-{Asymptotic} {Analysis} of {Robust} {Control} from
  {Coarse}-{Grained} {Identification}.
\newblock \emph{arXiv:1707.04791 [cs, math]}.
\newblock \urlprefix\url{http://arxiv.org/abs/1707.04791}.
\newblock ArXiv: 1707.04791.

\bibitem[{Wagenmaker and Jamieson(2020)}]{wagenmaker_active_2020}
Wagenmaker, A. and Jamieson, K. (2020).
\newblock Active {Learning} for {Identification} of {Linear} {Dynamical}
  {Systems}.
\newblock In \emph{Proceedings of {Thirty} {Third} {Conference} on {Learning}
  {Theory}}, 3487--3582. PMLR.
\newblock
  \urlprefix\url{https://proceedings.mlr.press/v125/wagenmaker\\20a.html}.
\newblock ISSN: 2640-3498.

\bibitem[{Yin et~al.(2021)Yin, Iannelli, and Smith}]{yin_maximum_2021}
Yin, M., Iannelli, A., and Smith, R.S. (2021).
\newblock Maximum {Likelihood} {Estimation} in {Data}-{Driven} {Modeling} and
  {Control}.
\newblock \emph{IEEE Transactions on Automatic Control}, 1--1.
\newblock \doi{10.1109/TAC.2021.3137788}.
\newblock Conference Name: IEEE Transactions on Automatic Control.

\bibitem[{Zheng and Li(2021)}]{zheng_non-asymptotic_2021}
Zheng, Y. and Li, N. (2021).
\newblock Non-{Asymptotic} {Identification} of {Linear} {Dynamical} {Systems}
  {Using} {Multiple} {Trajectories}.
\newblock \emph{IEEE Control Systems Letters}, 5(5), 1693--1698.
\newblock \doi{10.1109/LCSYS.2020.3042924}.
\newblock \urlprefix\url{https://ieeexplore.ieee.org/document/9284539/}.

\bibitem[{Şen and Kalyoncu(2019)}]{gujs461494}
Şen, M.A. and Kalyoncu, M. (2019).
\newblock Grey wolf optimizer based tuning of a hybrid lqr-pid controller for
  foot trajectory control of a quadruped robot.
\newblock \emph{Gazi University Journal of Science}, 32(2), 674 -- 684.

\end{thebibliography}
                                                   






\appendix
\section{Proof of theorem~\ref{thm:asymp_error}}\label{append:asymp_error}
   
\begin{pf}
   Using the specific expression of $H_{k}$ in~\eqref{eq:Markov_parameter_form}, the bias of the estimation $\hat H_k$ has the following relationship with the polynomial approximation error:
   \begin{equation}\label{eq:asymp_error_decomp}
      \begin{aligned}
         \left|H_{k}-\sum_{t=1}^T\alpha_{t-1}^*H_t\right|&= \left|\sum_{j=1}^n c_j\lambda_j^{k-1}-\sum_{t=0}^{T-1} \alpha_{t}^*\sum_{j=1}^n c_j\lambda_j^t\right| \\
         & =\left|\sum_{j=1}^n c_j\left(\lambda_j^{k-1}-\sum_{t=0}^{T-1} \alpha_{t}^*\lambda_j^t\right)\right| \\
         & \leq C_m\left\|x^{k-1}-\sum_{t=0}^{T-1} \alpha_{t}^*x^t\right\|_{\infty, \rho},
      \end{aligned}
   \end{equation}
   where $\boldsymbol{\alpha}^*$ is obtained by~\eqref{eq:unregularized_uniform_approx}.
   \begin{itemize}
   \item Proof of $|H_k-\hat H_k|^2\leq 4C_m^2\rho^{2k-2}\exp\left(-\frac{2(T-1)^2}{k-1}\right)$
   \end{itemize}
   According to \cite{saibaba2021approximating}, the Chebyshev approximation problem has the following upper bound on the approximation error, where $\mu\in[-1, 1]$:
   \begin{equation}\label{eq:chebyshev_upper_bound}
      \left\|\mu^{k-1}-\sum_{t=0}^{T-1} \beta_{t}^*\mu^t\right\|_{\infty, 1}\leq 2\exp\left(-\frac{(T-1)^2}{k-1}\right),
   \end{equation}
   and the optimal solution $\beta_t^*$ is unique.
   Since $\mu\in[-1, 1]$, take $x=\rho \mu$, then $x\in[-\rho, \rho]$, and the following inequality holds by putting $\mu=\rho^{-1}x$ into~\eqref{eq:chebyshev_upper_bound}:
   \begin{equation}
      \begin{aligned}
         \left\|\rho^{1-k}x^{k-1}-\sum_{t=0}^{T-1}\rho^{-t}\beta^*_{t}x^t\right\|_{\infty, \rho} \leq 2\exp\left(-\frac{(T-1)^2}{k-1}\right).
      \end{aligned}
   \end{equation}
   Therefore, 
   \begin{equation}\label{eq:chebyshev_upper_bound_rho}
      \begin{aligned}
         \left\|x^{k-1}-\sum_{t=0}^{T-1}\rho^{k-t-1}\beta^*_tx^t\right\|_{\infty, \rho} \leq 2\rho^{k-1}\exp\left(-\frac{(T-1)^2}{k-1}\right).
      \end{aligned}
   \end{equation}
 Thus, the upper bound can be proved by combining~\eqref{eq:chebyshev_upper_bound_rho} with~\eqref{eq:asymp_error_decomp}.
   \begin{itemize}
      \item Proof of $|H_k-\hat H_k|^2\leq\frac{1}{4}C_m^2\rho^{2k-2}$
   \end{itemize}
   The bound can be derived by considering the following two cases. When $k$ is odd, 
   it can be verified that the following inequality holds:
   \begin{equation}
      \left|x^{k-1}-\frac{1}{2}\rho^{k-1}\right|\leq \frac{1}{2}\rho^{k-1}, \forall x\in[0, \rho].
   \end{equation}
   Since LHS is an even function, the polynomial approximation error can be upper bounded by:
   \begin{equation}\label{eq:even_upper_bound}
      \left\|x^{k-1}-\sum_{t=0}^{T-1} \alpha_{t}^*x^t\right\|_{\infty, \rho}\leq \left\|x^{k-1}-\frac{1}{2}\rho^{k-1}\right\|_{\infty, \rho}\leq \frac{1}{2}\rho^{k-1}.
   \end{equation}
   Similarly, when $k$ is even, the following inequality holds:
   \begin{equation}
      \left|x^{k-1}-\frac{1}{2}\rho^{k-2}x\right|\leq \frac{1}{2}\rho^{k-1}, \forall x\in[0, \rho].
   \end{equation}

   Thus,
   \begin{equation}\label{eq:odd_upper_bound}
      \left\|x^{k-1}-\sum_{t=0}^{T-1}\alpha_t^*x^t\right\|_{\infty, \rho}\leq \left\|x^{k-1}-\frac{1}{2}\rho^{k-2}x\right\|_{\infty, \rho}\leq \frac{1}{2}\rho^{k-1}.
   \end{equation}
   Therefore, the second upper bound in~\eqref{eq:asymp_error} can be derived by combining~\eqref{eq:asymp_error_decomp},~\eqref{eq:even_upper_bound} and~\eqref{eq:odd_upper_bound}.

\end{pf}

\section{Proof of Theorem~\ref{thm:markov_approx_error}}\label{append:markov_error}

\begin{pf}
   The mean squared approximation error could be expanded as:
   \begin{equation}
      \begin{aligned}
         &\mathbb{E}\left(H_k-\hat H_k\right)^2=\mathbb{E}\left(H_k-\sum_{t=1}^T\alpha_{t-1}\tilde H_t^{(N)}\right)^2 \\
         &=\mathbb{E}\left(H_k-\sum_{t=1}^T\alpha_{t-1}H_t+\sum_{t=1}^T \alpha_{t-1}H_t-\sum_{t=1}^T\alpha_{t-1}\tilde H_t^{(N)}\right)^2\\
         &=\left(H_k-\sum_{t=1}^T\alpha_{t-1}H_t\right)^2+\mathbb{E}\left[\sum_{t=1}^T\alpha_{t-1}\left(H_t-\tilde H_t^{(N)}\right)\right]^2,
      \end{aligned}
   \end{equation}
   which is the bias-variance decomposition.

   Therein, the bias can be written as follows according to~\eqref{eq:asymp_error_decomp}:
   \begin{equation}\label{eq:noisy_error_bias}
      \left(H_k-\sum_{t=1}^T\alpha_{t-1}H_t\right)^2\leq C_m^2\left\|x^{k-1}-\sum_{t=0}^{T-1}\alpha_tx^t\right\|_{\infty, \rho}^2.
   \end{equation}
   On the other hand, for random variables $X$ and $Y$, the following inequalities hold:
   \[ \text{Var}(X+Y)=\text{Var}(X)+\text{Var}(Y)+2\text{Cov}(X, Y) \]
   \begin{equation}
      \begin{aligned}
         & \leq \text{Var}(X)+\text{Var}(Y)+2\sqrt{\text{Var}(X)\text{Var}(Y)} \\
         & \leq \left(\sqrt{\text{Var}(X)}+\sqrt{\text{Var}(Y)}\right)^2.
      \end{aligned}
   \end{equation}
   Thus, the variance term satisfies:
   \begin{equation}\label{eq:markov_parameter_variance}
      \begin{aligned}
           &\mathbb{E}\left[\sum_{t=1}^T\alpha_{t-1}\left(H_t-\tilde H_t^{(N)}\right)\right]^2
          =\text{Var}\left[\sum_{t=1}^T\alpha_{t-1}\tilde H_t^{({N})}\right] \\
           & \leq \left(\sum_{t=1}^T \sqrt{\text{Var}\left(\alpha_{t-1}\tilde H_t^{(N)}\right)}\right)^2. 
      \end{aligned}
   \end{equation}
   By expanding $\tilde H_t^{(N)}$ as the mean of the $t$-th system output from $N$ independent experiments, we can obtain the following result:
   \begin{equation}
      \begin{aligned}
         \text{Var}\left(\tilde H_t^{(N)}\right)&=\text{Var}\left(\frac{1}{N}\sum_{\ell=1}^N y_t^{(\ell)}\right) \\
         & = \frac{1}{N^2}\sum_{\ell=1}^N\text{Var}\left(y_t^{(\ell)}\right) =\frac{\Sigma}{N},
      \end{aligned}
   \end{equation}
   where $\Sigma$ denotes the variance of system output $y_t$ in each episode, which is the same for each $t$ due to the fact that the input in the experiment is determined. Furthermore, since each episode is reset to steady state, $\Sigma$ is equal to:
   \begin{equation}\label{eq:markov_param_variance}
      \Sigma= CPC^\top+R, \text{where}\ P=APA^\top+Q.
   \end{equation}
   Therefore, equation~\eqref{eq:markov_parameter_variance} can be further derived as:
   \begin{equation}
      \left(\sum_{t=1}^T \sqrt{\alpha_{t-1}^2\text{Var}\left(\tilde H_t^{(N)}\right)}\right)^2  \leq \frac{\Sigma}{N}\left(\sum_{t=0}^{T-1}| \alpha_{t}|\right)^2 
          =\frac{\Sigma}{N}\|\boldsymbol{\alpha}\|_1^2.
   \end{equation}
   Thus, the theorem can be proved by combining~\eqref{eq:noisy_error_bias} with the result above.
\end{pf}

\end{document}